\documentclass[11pt, reqno]{amsart}
\usepackage{amsmath, amsthm, amscd, amsfonts, amssymb, graphicx, color, stmaryrd}
\usepackage[all,cmtip]{xy}	% diagrams
\usepackage[bookmarksnumbered, colorlinks, plainpages]{hyperref}

\newcommand{\Nn}{\mathbb{N}}
\newcommand{\Rr}{\mathbb{R}}
\newcommand{\Cc}{\mathbb{C}}
\newcommand{\Zz}{\mathbb{Z}}

\newcommand{\F}{\mathcal{F}}
\newcommand{\K}{\mathcal{K}}	% compact Operators
\renewcommand{\H}{\mathcal{H}}	% Hilbert-space
\renewcommand{\L}{\mathcal{L}}	% bd lin Operators
\newcommand{\C}{\mathcal{C}}	% Calkin-Algebra
\newcommand{\E}{\mathcal{E}}	% C^{\ast} algebras
\newcommand{\A}{\mathcal{A}} % Algebra bising. Pseudos 
\newcommand{\Ad}{\mathrm{Ad}}	
\newcommand{\Symb}{\Gamma}	% symbol classes
\renewcommand{\S}{\Symb}

\newcommand{\csigma}{\mathrm{\overline{\sigma}}}	% complete symbol
\newcommand{\tsigma}{\mathrm{\tilde{\sigma}}}		% symbol action
\newcommand{\symb}{\mathrm{symb}}
\newcommand{\tsymb}{\widetilde{\symb}}		% pt-wise Symbol

\newcommand{\Tau}{\mathcal{T}}
\newcommand{\id}{\mathrm{id}}
\newcommand{\ch}{\mathrm{ch}}			% Chern-character
\newcommand{\Ext}{\mathrm{Ext}}			% extensions

\newcommand{\scal}[2]{\langle #1, #2 \rangle}	% inner prod
\newcommand{\im}{\operatorname{im}} 
	
\newcommand{\coker}{\operatorname{coker}}	% Co-Kern
\newcommand{\op}{\operatorname{op}} 		% symbol-operator
\newcommand{\ind}{\operatorname{ind}}
\newcommand{\potimes}{\hat{\otimes}} % proj. tensor product

\newcommand{\iso}{\xrightarrow{\sim}}	% iso arrow

\newtheorem{Thm}{Theorem}[section]
\newtheorem{Lem}[Thm]{Lemma}
\newtheorem{Prop}[Thm]{Proposition}
\newtheorem{Cor}[Thm]{Corollary}
\theoremstyle{definition}
\newtheorem{Def}[Thm]{Definition}

\newtheorem*{acknowledgement}{Acknowledgement}
\newtheorem{Rem}[Thm]{Remark}

\newtheorem{Not}[Thm]{Notation}
\begin{document}
\setcounter{page}{1}

%-------------------------- Pleased do not change the following line-------------------------------------------
%\noindent \textcolor[rgb]{0.99,0.00,0.00}{This is a submission to one of journals of TMRG: BJMA/AFA}\\[.5in]
%--------------------------------------------------------------------------------------------------------------

\title{The $K$-theory of bisingular pseudodifferential operators}

\author[Karsten Bohlen]{Karsten Bohlen}

\address{$^{1}$ Leibniz University Hannover, Germany}
\email{\textcolor[rgb]{0.00,0.00,0.84}{bohlen.karsten@math.uni-hannover.de}}

%\dedicatory{This paper is dedicated to Professor ABCD}



% \date{Received: xxxxxx; Revised: yyyyyy; Accepted: zzzzzz.}

%\begin{abstract}
%\end{abstract} \maketitle

\begin{abstract}
This work is concerned with a class of pseudodifferential operators of tensor product type called bisingular operators. 
For the \emph{global} bisingular calculus in the flat ($\Rr^{n_1 + n_2}$) case we calculate the $K$-theory of $C^{\ast}$-algebras given by the norm-closures of spaces of bisingular pseudodifferential operators.
\end{abstract} 

\maketitle

\tableofcontents

% preliminary sections:

\section{Introduction}

The bisingular pseudodifferential operators were introduced in 1975 by L. Rodino. The bisingular calculus
defines pseudodifferential operators on the cartesian product $X \times Y$ of two compact smooth manifolds $X$ and $Y$.
This includes for example the tensor product $P_1 \otimes P_2$ as well as the external product $P_1 \sharp P_2$ of two classical
pseudodifferential operators. Together with the multiplicative property of the Fredholm index
the latter example was a motivation for the introduction of a general bisingular calculus. See e.g. \cite{rodino} and the references contained therein.
Later, a global bisingular calculus was introduced in the paper \cite{bgpr}, defined on the product $\Rr^{n_1} \times \Rr^{n_2}$ 
and extending the Shubin calculus.

In this paper we will continue the study of the bisingular pseudodifferential operators. 
The first goal is to determine the $K$-theory of the $C^{\ast}$-algebra completions of the global bisingular operators.
We thus consider the $C^{\ast}$-algebras of completed pseudodifferential operators obtained from the 
global bisingular calculus. Our main aim is to then calculate the $K$-theory of these algebras. 

For the case of the Shubin calculus of pseudodifferential operators defined on $\Rr^n$, the $K$-theory 
is easily calculated from an exact sequence which is induced by the principal symbol map. The bisingular calculus by contrast has two operator-valued symbols which take 
values in a non-commutative symbol algebra.
We will now explain more precisely the difference between these two cases.

First recall the construction of the standard exact sequence of pseudodifferential operators defined on $\Rr^n$.
Let $G_{cl}^m(\Rr^n)$ denote the classical pseudodifferential operators (operators of Shubin type) of order $m \in \Rr$.
Denote by $A$ the completion of $G_{cl}^0(\Rr^n)$ in the induced $\L(L^2)$ norm. We obtain a $C^{\ast}$ algebra $A := \overline{G_{cl}^{0}(\Rr^{n})}$ and 
the completion of $G_{cl}^{-1}$ yields $\K_{\Rr^n} = \K(L^2)$ the algebra of compact operators on $L^2(\Rr^n)$.
Denote by $\sigma \colon G_{cl}^0(\Rr^n) \to C^{\infty}(S^{2n-1})$ the principal symbol map and
$\overline{\sigma} \colon A \to C(S^{2n-1})$ its continuous extension. 
We have the following short exact sequence
\begin{align}
\xymatrix{
0 \ar[r] & \K_{\Rr^n} \ar[r] & A \ar[r]^-{\csigma} & C(S^{2n-1}) \ar[r] & 0.
} \label{S1}
\end{align} 

By comparison consider now the \emph{global bisingular} calculus.
Here we will consider algebras $G_{cl}^{m_1, m_2}(\Rr^{n_1} \times \Rr^{n_2})$ of the orders $m_1 \in \Rr$ and $m_2 \in \Rr$ which define
a global pseudodifferential calculus on products $\Rr^{n_1} \times \Rr^{n_2}$.
It turns out that taking completions with regard to $\L(L^2(\Rr^{n_1} \times \Rr^{n_2}))$ we find 
$\overline{G_{cl}^{-1,-1}(\Rr^{n_1} \times \Rr^{n_2})} = \K$ the compact operators on $L^2$.
If we set $\A := \overline{G_{cl}^{0,0}}$ we will prove the exactness of the sequence
\begin{align}
\xymatrix{
0 \ar[r] & \K \ar[r] & \A \ar[r]^-{\overline{\sigma}_1 \oplus \overline{\sigma}_2} & \Sigma \ar[r] & 0 
} \label{S2}
\end{align}

with a non-commutative symbol-algebra $\Sigma$ and the extended direct sum principal symbol map $\overline{\sigma}_1 \oplus \overline{\sigma}_2$.

The algebra $\Sigma$ can be viewed as a $C^{\ast}$-algebra pullback, i.e. a restricted direct sum of 

$C(S^{2n_1 -1}, \overline{G_{cl}^0(\Rr^{n_2})})$ and $C(S^{2n_2-1}, \overline{G_{cl}^0(\Rr^{n_1})})$.
We will use this understanding of the symbol algebra to derive its $K$-theory. 

Our first main result can thus be stated:
\begin{align*}
K_0(\Sigma) &\cong \ker(\tsigma_{{\Rr^{n_1}}^{\ast}} - \tsigma_{{\Rr^{n_2}}^{\ast}}) \cong \Zz, \\
K_1(\Sigma) &\cong \coker(\tsigma_{{\Rr^{n_1}}^{\ast}} - \tsigma_{{\Rr^{n_2}}^{\ast}}) \cong \Zz.
\end{align*}

With isomorphisms induced by continuous extensions of pointwise quotient mappings (denoted here by $\tsigma$).

Another major theme which is not well-understood in the bisingular case is index theory. 
We first analyse the index problem for global bisingular operators. 
Then we pose a general index problem which does not depend on the calculus. 
We consider the index theory of tensor products of three classes of Toeplitz operators.
An interesting observation is that the index theory for bisingular operators is the same as for these classes of Toeplitz operators.

\begin{acknowledgement}
I thank Elmar Schrohe for many helpful discussions and for bringing these problems to my attention.
Secondly, I'm grateful to Magnus Goffeng for numerous discussions on index theory.
Part of this work was conducted while I was a member of the GRK 1463 at Leibniz University of Hannover.
I thank the Deutsche Forschungsgemeinschaft (DFG) for their financial support.
\end{acknowledgement}

\section{Bisingular operators}

% the class for compact smooth manifolds X_1, X_2
% all cool properties...
% relation to classical and non-classical symbols

In this section we will introduce the terminology and notation for the rest of the paper.

\begin{Def}
The class of bisingular symbols $\Gamma^{m_1, m_2}(\Rr^{n_1 + n_2}; \Cc^{l \times k})$ consists of smooth functions
$a \colon \Rr^{2n_1 + 2n_2} \to \Cc^{l \times k}$ with the following uniform estimates
\[
\|D_{\xi_1}^{\alpha_1} D_{x_1}^{\beta_1} D_{\xi_2}^{\alpha_2} D_{x_2}^{\beta_2} a(x_1, \xi_1, x_2, \xi_2)\|_{\Cc^{l \times k}} \leq C \scal{x_1}{\xi_1}^{m_1 - |\alpha_1| - |\beta_1|} \scal{x_2}{\xi_2}^{m_2 - |\alpha_2| - |\beta_2|}. 
\]

Furthermore, we set
\[
\Gamma^{-\infty, -\infty}(\Rr^{n_1 + n_2}; \Cc^{l \times k}) = \bigcap_{m_1, m_2} \Gamma^{m_1, m_2}(\Rr^{n_1 + n_2}; \Cc^{l \times k}). 
\]

\label{Def:global}
\end{Def}

% introduce classical symbols

% tensor product
% algebraic structure
% symbol space definition

Given such a symbol $a$ we have two maps
\[
(x_1, \xi_1) \mapsto a_1(x_1, \xi_1) := ((x_2, \xi_2) \mapsto a(x_1, \xi_1, x_2, \xi_2)) 
\]

and
\[
(x_2, \xi_2) \mapsto a_2(x_2, \xi_2) := ((x_1, \xi_1) \mapsto a(x_1, \xi_1, x_2, \xi_2)). 
\]

Hence $a_1 \in \Gamma^{m_2}(\Rr^{n_2}, \Gamma^{m_1}(\Rr^{n_1}; \Cc^{l \times k}))$ 
and $a_2 \in \Gamma^{m_1}(\Rr^{n_1}, \Gamma^{m_2}(\Rr^{n_2}; \Cc^{l \times k}))$. 

The subclass of bisingular classical symbols is denoted by $\Gamma_{cl}^{m_1, m_2}$ and obtained by using in the above definition the classical Shubin classes.

% describe asymptotic expansion etc.

We have two principal symbols
\begin{align}
\sigma_1^{m_1}(A) &= a_1^{(m_1)} \in C^{\infty}(S^{2n_1 - 1}, G_{cl}^{m_2}(\Rr^{n_2})), \label{p1} \\
\sigma_2^{m_2}(A) &= a_2^{(m_2)} \in C^{\infty}(S^{2n_2 - 1}, G_{cl}^{m_1}(\Rr^{n_1})). \label{p2} 
\end{align}

The principal symbols have the following properties for $A \in G_{cl}^{m_1, m_2}(\Rr^{n_1 + n_2}), \ B \in G_{cl}^{p_1, p_2}(\Rr^{n_1 + n_2})$

\begin{align*}
\sigma_i^{m_i + p_i}(A \cdot B) = \sigma_i^{m_i}(A) \cdot \sigma_i^{p_i}(B), \\ 
\sigma_i^{m_i}(A^{\ast}) = \sigma_i^{m_i}(A)^{\ast}, i = 1,2. 
\end{align*}

% joint symbol

% ptwise principal symbols

Fix the notation $\sigma_{\Rr^{n_1}}, \sigma_{\Rr^{n_2}}$ for the principal symbol map of $G_{cl}^{m_1}(\Rr^{n_1})$ and $G_{cl}^{m_2}(\Rr^{n_2})$
respectively. 
Then define in each case the pointwise principal symbol maps
\begin{align*}
&\tilde{\sigma}_{\Rr^{n_1}} \colon C^{\infty}(S^{2n_2 - 1}, G_{cl}^{m_1}(\Rr^{n_1})) \to C^{\infty}(S^{2n_1 - 1} \times S^{2n_2 - 1}), \\
&\tilde{\sigma}_{\Rr^{n_1}}(F)(x_1, \xi_1, x_2, \xi_2) := \sigma_{\Rr^{n_1}}(F(x_2, \xi_2))(x_1, \xi_1), \ F \in C^{\infty}(S^{2n_2-1}, G_{cl}^{m_1}(\Rr^{n_1})), \\
&\tilde{\sigma}_{\Rr^{n_2}} \colon C^{\infty}(S^{2n_1 - 1}, G_{cl}^{m_2}(\Rr^{n_2})) \to C^{\infty}(S^{2n_1 - 1} \times S^{2n_2 - 1}), \\
&\tilde{\sigma}_{\Rr^{n_1}}(G) := \sigma_{\Rr^{n_2}}(G(x_1, \xi_1))(x_2, \xi_2), \ G \in C^{\infty}(S^{2n_1 - 1}, G_{cl}^{m_2}(\Rr^{n_2}).  
\end{align*}

Note that by nuclearity we have 
\begin{align*}
&C^{\infty}(S^{2n_1 - 1}, G_{cl}^{m_2}(\Rr^{n_2})) \cong C^{\infty}(S^{2n_1 - 1}) \potimes G_{cl}^{m_2}(\Rr^{n_2}), \\
&C^{\infty}(S^{2n_2 - 1}, G_{cl}^{m_1}(\Rr^{n_1})) \cong C^{\infty}(S^{2n_2 - 1}) \potimes G_{cl}^{m_1}(\Rr^{n_1}) 
\end{align*}

and the pointwise symbol maps are also given by
\[
\tilde{\sigma}_{\Rr^{n_1}} = \id_{C^{\infty}(S^{2n_2 -1})} \otimes \sigma_{\Rr^{n_1}}, \ \tilde{\sigma}_{\Rr^{n_2}} = \sigma_{\Rr^{n_2}} \otimes \id_{C^{\infty}(S^{2n_1 - 1})}. 
\]

We assume the \emph{compatibility condition}
\begin{align}
& \sigma_{\Rr^{n_2}}(\sigma_1^{m_1}(A)(x_1, \xi_1))(x_2, \xi_2) = \sigma_{\Rr^{n_1}}(\sigma_2^{m_2}(A)(x_2, \xi_2))(x_1, \xi_1) \notag \\
&= \sigma^{m_1, m_2}(A)(x_1, \xi_1, x_2, \xi_2) = a_{m_1, m_2}(x_1, \xi_1, x_2, \xi_2). \label{comp}
\end{align}

\begin{Def}
Let $\Sigma^{m_1, m_2}$ be the set of all pairs 
\[
(F, G) \in C^{\infty}(S^{2n_1-1}, G_{cl}^{m_1}(\Rr^{n_2})) \oplus C^{\infty}(S^{2n_2-1}, G_{cl}^{m_2}(\Rr^{n_1})) 
\]
such that 
\begin{align*}
\tsigma_{\Rr^{n_2}}(F) &= \tsigma_{\Rr^{n_1}}(G).
\end{align*}

Let $\{F_1, G_1\} \in \Sigma^{m_1, m_2}, \ \{F_2, G_2\} \in \Sigma^{p_1, p_2}$ and set
\begin{align*}
\{F_2, G_2\} \circ \{F_1, G_1\} &:= \{F_2 \circ_{2} F_1, G_2 \circ_{1} G_1\} \in \Sigma^{m_1 + p_1, m_2 + p_2}.
\end{align*}

Here
\begin{align*}
(F_2 \circ_{2} F_1)(x_1, \xi_1) &:= F_2(x_1, \xi_1) \circ_{\Rr^{n_2}} F_1(x_1, \xi_1), \\
(G_2 \circ_1 G_1)(x_2, \xi_2) &:= G_2(x_2, \xi_2) \circ_{\Rr^{n_1}} G_1(x_2, \xi_2)
\end{align*}

where $\circ_{\Rr^{n_2}}$ denotes the operator product $G_{cl}^{m_2}(\Rr^{n_2}) \times G_{cl}^{p_2}(\Rr^{n_2}) \to G_{cl}^{m_2 + p_2}(\Rr^{n_2})$ and for $\circ_{\Rr^{n_1}}$ analogously.
\label{Def:symbspace}
\end{Def}

\begin{Prop}
The following sequence is exact
\begin{align}
\xymatrix{
G_{cl}^{m_1 - 1,m_2 - 1}(\Rr^{n_1 + n_2}) \ar@{>->}[r]^-{i} & G_{cl}^{m_1, m_2}(\Rr^{n_1 + n_2}) \ar@{->>}[r]^-{\sigma_1^{m_1} \oplus \sigma_2^{m_2}} & \Sigma^{m_1, m_2}. 
} \label{exa}
\end{align}

\label{Prop:exa}
\end{Prop}

\begin{proof}
\emph{1.} Let $P = \op(a), \ a \in \Gamma_{cl}^{m_1, m_2}$ such that $\sigma_1^{m_1}(P) = 0, \ \sigma_2^{m_2}(P) = 0$.
Then it follows from Def. 1.6. \cite{bgpr} (cf. \cite{rodino}, Def. 2.1.) that $a \in \Gamma_{cl}^{m_1 - 1, m_2}$ and also $a \in \Gamma_{cl}^{m_1, m_2 - 1}$. 
Now we use that $a$ is classical in the sense that we can find $a_k$ homogenous of degree $m_1 - k$ in $(x_1, \xi_1) \in \Rr^{2n_1}$ for $k = 0, \cdots, N$ such that
\begin{align*}
& a - \sum_{k=0}^N a_k \in \Gamma^{m_1 - (N+1),m_2}
\end{align*}

for each $N \in \Nn$. And $b_k$ homogenous of degree $m_2 - k$ for $k = 0, \cdots, N$ in $(x_2, \xi_2) \in \Rr^{2n_2}$ such that
\begin{align*}
& a - \sum_{k=0}^N b_k \in \Gamma^{m_1, m_2 - (N+1)}
\end{align*}

for each $N \in \Nn$. 

Therefore by considering $k = 0$ we can write
\begin{align*}
a &= \tilde{b}_0 + \tilde{b} = \tilde{a}_0 + \tilde{a}
\end{align*}

for $\tilde{a}, \ \tilde{b} \in \Gamma^{m_1 - 1, m_2 - 1}$ as well as $\tilde{b}_0 \in \Gamma^{m_1 - 1, m_2}$ homogenous in $(x_2, \xi_2)$ of order $m_2$ and
$\tilde{a}_0 \in \Gamma^{m_1, m_2 - 1}$ homogenous in $(x_1, \xi_1)$ of order $m_1$. 
Then $\tilde{a}_0 = \tilde{b}_0 + \tilde{b} - \tilde{a} \in \Gamma^{m_1 - 1, m_2 - 1}$ and hence $\tilde{a}_0$ is also of order $m_1 - 1$.
But this implies that $\tilde{a}_0 = 0$ and analogously, $\tilde{b}_0 = 0$. 
It follows that $a \in \Gamma_{cl}^{m_1 - 1, m_2 -1}$. 

Now for $P = \op(a), \ a \in \Gamma_{cl}^{m_1 - 1, m_2 - 1}$ we must have $\sigma_1(P) = 0, \ \sigma_2(P) = 0$ as 
$\Gamma_{cl}^{m_1 - 1, m_2 - 1} \subseteq \Gamma_{cl}^{m_1 -1, m_2} \cap \Gamma_{cl}^{m_1, m_2 - 1}$.
It follows that $\ker(\sigma_1 \oplus \sigma_2) = \ker \sigma_1 \cap \ker \sigma_2 = \im(i)$.

\emph{2.} Let $(F, G) \in \Gamma_{cl}^{m_1, m_2}$. There are two maps
\[
\symb_i \colon G_{cl}^{m_i}(\Rr^{n_i}) \to \Gamma_{cl}^{m_i}(\Rr^{n_i}), \ i = 1,2
\]

which are right-inverse to $\op_i \colon \Gamma_{cl}^{m_i}(\Rr^{n_i}) \to G_{cl}^{m_i}(\Rr^{n_i}), \ i = 1,2$ modulo smoothing terms (defined by
$A \mapsto e^{-ix_i \xi_i} A e^{i x_i \xi_i}$).
Set $r := \tsigma_{\Rr^{n_2}}(F) \in C^{\infty}(S^{2n_1-1} \times S^{2n_2-1})$ and 
\[
p := \tsymb_2 \circ F, \ q := \tsymb_1 \circ G.
\]

Where $\tsymb_i$ are the pointwise evaluations of $\symb_i$ for $i = 1,2$.
Choose two smooth cut-off functions near $0$ on $\Rr^{2n_1}, \Rr^{2n_2}$ respectively: $\chi_1, \ \chi_2$ and set
\[
a := \chi_1 p + \chi_2 q - \chi_1 \chi_2 r \in \S^{m_1, m_2}.
\] 

Then $\sigma_1(\op(a)) = F$ and $\sigma_2(\op(a)) = G$ (cf. \cite{bgpr}, Def. 1.6 iii)). 
Hence $\sigma_1 \oplus \sigma_2$ is surjective and the exactness of \eqref{exa} is established. 
\end{proof}

We introduce the appropriate Sobolev spaces for bisingular operators.
\begin{Def}
We define the Sobolev space as the completion 
\[
Q^{s,t}(\Rr^{n_1 + n_2}) = \overline{S(\Rr^{2n_1 + 2n_2})}^{\|\cdot\|_{s,t}}
\]

where the norm is given by
\[
\|u\|_{s,t} := \|\Lambda^{s,t} u\|_{L^2(\Rr^{n_1 + n_2})}. 
\]

Here $\Lambda^{s,t} := \Lambda_{n_1}^{s} \otimes \Lambda_{n_2}^{t}$ and $\Lambda_{n_i}$ are invertible operators
on $L^2(\Rr^{n_1}), \ i=1,2$. 
\label{Def:Sobolev}
\end{Def}

\begin{Prop}[cf. \cite{bgpr}]
Let $P \in G_{cl}^{m_1, m_2}(\Rr^{n_1 + n_2}, \Cc^{l \times k})$ then $P$ has a continuous linear extension $P \colon Q^{s, t}(\Rr^{n_1+n_2}, \Cc^{k}) \to Q^{s-m_1, t-m_2}(\Rr^{n_1 + n_2}, \Cc^{l})$. 
\label{Prop:cont}
\end{Prop}

% continuity on Sobolev space

\section{$K$-theory}
\label{Kthy}

% introduce globally bisingular operators
% completions, top tensor product for the classical
% all the rest.. including pullback MV calculation
In what follows we want to consider the norm completions of the classical, global bisingular operators of orders $(0, 0)$ as previously introduced.

\begin{Lem}
\emph{i)} Completion with respect to the bounded linear operators on $L^2(\Rr^{n_1} \times \Rr^{n_2})$ yields an isomorphism
\begin{align}
\overline{G_{cl}^{-1, -1}(\Rr^{n_1} \times \Rr^{n_2})}^{L^2} &\cong \K(L^2(\Rr^{n_1} \times \Rr^{n_2})) \label{C1} 
\end{align}

\emph{ii)} We have the isomorphism 
\[
G_{cl}^{m_1}(\Rr^{n_1}) \potimes G_{cl}^{m_2}(\Rr^{n_2}) \cong G_{cl}^{m_1, m_2}(\Rr^{n_1} \times \Rr^{n_2}). 
\]

\label{Lem:1}
\end{Lem}

\begin{proof}
\emph{i)} We repeat the standard argument for the benefit of the reader.

We have $G^{-\infty, -\infty}(\Rr^{n_1 + n_2}) \subset G^{-1,-1}(\Rr^{n_1 + n_2})$ and operators in $G^{-\infty, -\infty}$
have integral kernel in $S(\Rr^{2n_1 + 2n_2})$. 
Now $S(\Rr^{2n_1 + 2n_2})$ is dense in $L^2(\Rr^{2n_1 + 2n_2})$.
Thus $G^{-\infty,-\infty}$ is dense in the Hilbert-Schmidt operators $\L^2(L^2(\Rr^{n_1 + n_2}))$. 
These are dense in the ideal of compact operators.

\emph{ii)} First of all we have a dense inclusion $\Gamma_{cl}^{m_1}(\Rr^{n_1}) \otimes \Gamma_{cl}^{m_2}(\Rr^{n_2}) \hookrightarrow \Gamma_{cl}^{m_1, m_2}(\Rr^{n_1 + n_2})$.
To see this exhibit a dense subspace on both sides, e.g. the homogenous smooth functions. 
By the nuclearity of the spaces the projective tensor product is the unique completed topological tensor product.
Denote the projective topology by $\pi$ as usual and by $\epsilon$ the injective topology. 
Now $(a, b) \mapsto a \otimes b$ is separately continuous as a map $\Gamma_{cl}^{m_1} \otimes \Gamma_{cl}^{m_2} \to \Gamma_{cl}^{m_1 + m_2}$.
From this we obtain that $\Gamma_{cl}^{m_1} \otimes \Gamma_{cl}^{m_2}$ induces on $\Gamma_{cl}^{m_1, m_2}$ a topology 
which is weaker than $\pi = \epsilon$. 
Then using the density it is a standard argument to prove that also $\Gamma_{cl}^{m_1} \otimes \Gamma_{cl}^{m_2}$ induces on 
$\Gamma_{cl}^{m_1, m_2}$ a stronger topology than $\pi = \epsilon$, cf. \cite{treves}. 
\end{proof}

\begin{Not}
Throughout the rest of this text we fix the following notation.

Introduce the respective $\L(L^2)$ completions
\begin{align*}
A_1 &:= \overline{G_{cl}^0(\Rr^{n_1})}, \ A_2 := \overline{G_{cl}^0(\Rr^{n_2})}, \ \K_1 := \overline{G_{cl}^{-1}(\Rr^{n_1})}, \ \K_2 := \overline{G_{cl}^{-1}(\Rr^{n_2})}, \\
A^{i,j} &:= \overline{G_{cl}^{i,j}(\Rr^{n_1} \times \Rr^{n_2})}, \ i, j \in \{-1,0,1\}.
\end{align*}

For $i = j = -1$ by the Lemma we have $A^{-1,-1} = \K$, the compact operators on the Hilbert space $L^2(\Rr^{n_1} \times \Rr^{n_2})$.
\label{Not:completions}
\end{Not}

\begin{Lem}
Introduce the $C^{\ast}$ subalgebra $\Sigma \subset C(S^{2n_1 - 1}, A_2) \oplus C(S^{2n_2 - 1}, A_1)$ 
\[
\Sigma := \{(F, G) : q_1(F) = q_2(G) \in C(S^{2n_1 - 1} \times S^{2n_2 - 1})\}.
\]

We set here $q_i$ for the canonical, pointwise quotient maps.
Then the $C^{\ast}$ completion $\overline{\Sigma^{0,0}}$ is isomorphic to $\Sigma$.
\label{Lem:compl}
\end{Lem}

\begin{proof}
We have the completions $\overline{C^{\infty}(S^{2n_1 - 1} \times S^{2n_2 - 1})} = C(S^{2n_1 - 1} \times S^{2n_2 - 1})$ and
$\overline{C^{\infty}(S^{2n_1 - 1}, G^0(\Rr^{n_2}) } = C(S^{2n_1 - 1}, A_2), \ \overline{C^{\infty}(S^{2n_2 - 1}, G^0(\Rr^{n_1}))} = C(S^{2n_2 - 1}, A_1)$.

Let $(F, G) \in \Sigma$ and choose a sequence $(f_n)_{n \in \Nn}$ in $C^{\infty}(S^{2n_1 - 1} \times S^{2n_2 - 1})$ such that
\begin{align*}
& \|q_2(G) - f_n\|_{\infty} = \|q_1(F) - f_n\|_{\infty} \leq \frac{1}{n}, \ n \in \Nn.
\end{align*}

Using the quotient topology with these estimates, density and surjectivity of quotient maps we can find sequences $(F_n)$ and $(G_n)$ in $C^{\infty}(S^{2n_1 - 1}, G^0(\Rr^{n_2}))$ and
$C^{\infty}(S^{2n_2 - 1}, G^0(\Rr^{n_1}))$ respectively such that
\begin{align*}
& \|F - F_n\|_{C(S^{2n_1 - 1}, A_2)} \leq \frac{1}{n}, \ n \in \Nn, \\
& \|G - G_n\|_{C(S^{2n_2 - 1}, A_1)} \leq \frac{1}{n}, \ n \in \Nn
\end{align*}

and such that
\begin{align*}
q_1(F - F_n) &= q_1(F - I_1\cdot f_n) = q_1(F) - f_n, \\
q_2(G - G_n) &= q_2(G - I_2 \cdot f_n) = q_2(G) - f_n 
\end{align*}

hence
\begin{align*}
q_1(F_n) = q_2(G_n) = f_n. 
\end{align*}

Since $q_1, \ q_2$ are continuous extensions we have that $(F_n, \ G_n) \in \Sigma^{0,0}$ for each $n \in \Nn$. 
\end{proof}

\begin{Cor}
The completion of $A_1 \otimes A_2$ with respect to any $C^{\ast}$-tensor norm is isomorphic to $\A$.
Also the completion of $\K_1 \otimes \K_2$ is isomorphic to $\K$.
\label{Cor:tensor}
\end{Cor}

\begin{proof}
We endow the algebraic tensor product $A_1 \otimes A_2$ with the spatial tensor product norm $\sigma$ (c.f. \cite{wegge}, Def. T.5.16). Since we have two injective $\ast$-representations
$A_1 \to \L(L^2(\Rr^{n_1})), A_2 \to \L(L^2(\Rr^{n_2}))$ given by the inclusions it follows that the dense homomorphism
$A_1 \otimes A_2 \to \A$ (by \ref{Lem:1}, \emph{ii)}) is an isometry from $(A_1 \otimes A_2, \sigma) \to (\A, \|\cdot\|)$.
Thus $A_1 \otimes_{\sigma} A_2 \cong \A$. 
As extensions of nuclear $C^{\ast}$-algebras $A_1$ and $A_2$ are nuclear (see e.g. \cite{wegge}, Thm. T.6.27). 
Therefore the isomorphism holds for any $C^{\ast}$ norm on the tensor product. 
The same argument applies to the case $\K_1 \otimes_{\sigma} \K_2 \cong \K$.
\end{proof}

From now on we write $A_1 \otimes A_2$ for the $C^{\ast}$-tensor product and identify 
$A_1 \otimes A_2$ with $\A$. 
In the same way we identify $\K_1 \otimes \K_2$ with $\K$. 

In view of the continuous extensions $\sigma_{\Rr^{n_1}} \colon A_1 \to C(S^{2n_1 - 1}), \ \sigma_{\Rr^{n_2}} \colon A_2 \to C(S^{2n_2 - 1})$ the pointwise principal symbol maps are given by
\begin{align*}
& \tsigma_{\Rr^{n_1}} \colon C(S^{2n_1 - 1}, A_2) \to C(S^{2n_1 - 1} \times S^{2n_2 - 1}), \\
& \tsigma_{\Rr^{n_2}} \colon C(S^{2n_2 - 1}, A_1) \to C(S^{2n_1 - 1} \times S^{2n_2 - 1}).
\end{align*}

The completed symbol-algebra $\Sigma$ is obtained as the restricted direct sum of $C(S^{2n_1 - 1}, A_2) \oplus C(S^{2n_2 - 1}, A_1)$ by Lemma \ref{Lem:compl}.
First we fix the projections $\pi_1, \pi_2$ from $\Sigma^{0,0}$ onto the first, respectively second component.
We denote the continuous extensions of these projections by the same symbols (for simplicity). 

With this $\Sigma$ is a $C^{\ast}$-algebra with norm:
\begin{align*}
\|(F, G)\| &:= \sup\{\|F\|_1, \|G\|_2\}, \ (F, G) \in \Sigma, \\
\|F\|_1 &:= \sup_{(x_2, \xi_2) \in S^{2n_2 - 1}} \|F(x_2, \xi_2)\|_{A_1}, \\
\|G\|_2 &:= \sup_{(x_1, \xi_1) \in S^{2n_1 - 1}} \|G(x_1, \xi_1)\|_{A_2}.
\end{align*}

The pullback $\Sigma$ can thus be written in terms of the following diagram:
\begin{align}
\xymatrix{
\Sigma \ar[d]_-{\pi_2} \ar[r]^-{\pi_1} & C(S^{2n_2 - 1}) \otimes A_1 \ar[d]_{\tsigma_{\Rr^{n_1}}} \\
C(S^{2n_1 - 1}) \otimes A_2 \ar[r]^{\tsigma_{\Rr^{n_2}}} & C(S^{2n_1 - 1} \times S^{2n_2 - 1})
} \label{PB1}
\end{align}

We obtain with this the following result.
\begin{Thm}
We have an isomorphism $\A / \K \cong \Sigma$ induced by the continuous extension $\sigma$ of the 
direct-sum principal symbol.
\label{Thm:exa2}
\end{Thm}

\begin{proof}
First consider the two short exact sequences:
\[
\xymatrix{
0 \ar[r] & \K_1 \otimes A_2 \ar[r] & \A \ar[r]^-{\sigma_{\Rr^{n_1}} \otimes \id} & C(S^{2n_1 - 1}) \otimes A_2 \ar[r] & 0, \\
0 \ar[r] & \K_2 \otimes A_1 \ar[r] & \A \ar[r]^-{\id \otimes \sigma_{\Rr^{n_2}}} & C(S^{2n_2 - 1}) \otimes A_1 \ar[r] & 0.
} 
\]

This is already enough information to construct the pullback $\Sigma$ (up to isomorphism) and identify it with $\A / \K$.

Consider the following diagram which is put together by tensoring of the standard exact sequences and application of quotient mappings:
\begin{align*}
\xymatrix{
& & & C(S^{2n_1 - 1}, \K_2) \ar@{>->}[rd] & \\
& & \K_2 \otimes A_1 \ar@{->>}[ur] \ar@{>->}[dr] & & C(S^{2n_1 - 1}, A_2) \ar@{->>}[rd]^{\tilde{\sigma}_{\Rr^{n_2}}} & \\
& \K_1 \otimes \K_2 \ar@{>->}[dr] \ar@{>->}[ur] \ar@{>->}[r(1.8)] & & A_1 \otimes A_2 \ar@{->>}[ur]^{\sigma_{\Rr^{n_1}} \otimes \id} \ar@{->>}[dr]_{\id \otimes \sigma_{\Rr^{n_2}}} \ar@{-->>}[r]^{q} & \widetilde{\Sigma} \ar@{-->}[d]_-{\tilde{\pi}_2} \ar@{-->}[u]_-{\tilde{\pi}_1} & C(S^{2n_1 - 1} \times S^{2n_2 - 1}) & \\
& & \K_1 \otimes A_2 \ar@{->>}[dr] \ar@{>->}[ur] & & C(S^{2n_2 - 1}, A_1) \ar@{->>}[ur]^{\tilde{\sigma}_{\Rr^{n_1}}} & \\
& & & C(S^{2n_2 - 1}, \K_1) \ar@{>->}[ur] & 
}
\end{align*}

Here we denote by $\widetilde{\Sigma}$ a pullback 
\[
C(S^{2n_1 - 1}, A_2) \oplus_{C(S^{2n_1 - 1} \times S^{2n_2 - 1})} C(S^{2n_2 - 1}, A_1).
\] 
The map $q$ is well-defined as follows: 
\[
q(x) := (\sigma_{\Rr^{n_1}} \otimes \id)(x) \oplus (\id \otimes \sigma_{\Rr^{n_2}})(x), \ x \in A_1 \otimes A_2
\]
with kernel 
\[
\ker q = \ker (\sigma_{\Rr^{n_1}} \otimes \id) \cap \ker (\sigma_{\Rr^{n_2}} \otimes \id) = (\K_1 \otimes A_2) \cap (A_1 \otimes \K_2) = \K_1 \otimes \K_2.
\]
Also $q$ is surjective: let $(F, G) \in \tilde{\Sigma}$ and choose $x \in A_1 \otimes A_2$ such that $(\sigma_{\Rr^{n_1}} \otimes \id)(x) = F$.
Then
\[
\tsigma_{\Rr^{n_2}}((\sigma_{\Rr^{n_1}} \otimes \id)(x) - G) = \tsigma_{\Rr^{n_1}}((\id \otimes \sigma_{\Rr^{n_2}})(x)) - \tsigma_{\Rr^{n_2}}(G) = \tsigma_{\Rr^{n_1}}(F) - \tsigma_{\Rr^{n_2}}(G) = 0
\]

which implies $(\sigma_{\Rr^{n_1}} \otimes \id)(x) - G \in \ker \tsigma_{\Rr^{n_2}} = C(S^{2n_1 -1}, \K_2)$.
But $C(S^{2n_1 - 1}, \K_2) = (\sigma_{\Rr^{n_1}} \otimes \id)(A_1 \otimes \K_2)$. 
Hence we find $x_0 \in A_1 \otimes \K_2$ with $(\sigma_{\Rr^{n_1}} \otimes \id)(x + x_0) = G$.  
It follows that $q(x + x_0) = (F, G)$.

Note that by uniqueness of the $C^{\ast}$ pullback $\Sigma$ is isomorphic to $\widetilde{\Sigma}$. 

Now the continuous extension of $\sigma_1 \oplus \sigma_2$ is denoted by $\sigma$. 

We see that $\sigma$ is surjective because $(\sigma_1 \oplus \sigma_2)(P^{\ast}) = (\sigma_1 \oplus \sigma_2)(P)^{\ast}$ 
and $\sigma_1 \oplus \sigma_2$ is surjective. Hence for a given $b \in \Sigma$ we find a sequence $(b_n)$ in $\Sigma^{0,0}$ converging in $C^{\ast}$ norm
to $b$. By surjectivity we find a sequence $(a_n)$ in $G_{cl}^{0,0}$ such that $\sigma_1 \oplus \sigma_2(a_n) = b_n$. By continuity of the extension
$(a_n)$ $C^{\ast}$ converges. So $a_n \to a$ und $a$ is contained in the closure, again by continuity $\sigma(a) = b$. 
Also $\K \subset \ker \sigma$ by Prop. \ref{Prop:exa}. 
Then $q$ restricted to $G_{cl}^{0,0}$ agrees with $\sigma_1 \oplus \sigma_2$ by the construction above which implies $\ker \sigma \subset \K$.
% hier vlt. ein kleines Diagram mit stet. Forts.
\end{proof}

\begin{Rem}
As the isomorphism induced by $\sigma_1 \oplus \sigma_2$ (continuous extensions) is automatically isometric it furnishes the norm-equality
\[
\inf_{K \in \K} \|P + K\| = \sup \{\|\sigma_1(P)\|_1, \|\sigma_2(P)\|_2\}, \ P \in \A
\]

which is an expected standard result for a pseudodifferential calculus.

\label{Rem:norm}
\end{Rem}

Next we will calculate the $K$-theory of the completed algebras.

\begin{Thm}
We have the following $K$-theory
\begin{align*}
K_0(A^{i,j}) &\cong \Zz, \ K_1(A^{i,j}) \cong 0, \ i, j = 0, -1 \\ 
K_0(\Sigma) &\cong \ker(\tsigma_{{\Rr^{n_1}}^{\ast}} - \tsigma_{{\Rr^{n_2}}^{\ast}}) \cong \Zz \\
K_1(\Sigma) &\cong \coker(\tsigma_{{\Rr^{n_1}}^{\ast}} - \tsigma_{{\Rr^{n_2}}^{\ast}}) \cong \Zz
\end{align*}

where we set $\tsigma_{{\Rr^{n_1}}^{\ast}} := K_0(\tsigma_{\Rr^{n_1}}), \ \tsigma_{{\Rr^{n_2}}^{\ast}} := K_0(\tsigma_{\Rr^{n_2}})$ for the induced maps in $K$-theory.

\label{Thm:Kthy}
\end{Thm}

\begin{proof}
\emph{i)} We first note that $K_0(A_1) = K_0(A_2) \cong \Zz, K_1(A_1) = K_1(A_2) \cong 0$. 
Using that $K_i(C(S^{2n_j-1})) = \Zz, i, j = 0,1$ (c.f. \cite{wegge}, 6.5) this follows by application of the six-term exact sequence applied to 
\begin{align*}
\xymatrix{
0 \ar[r] & \K_1 \ar[r] & A_1 \ar[r]^-{\sigma_{\Rr^{n_1}}} & C(S^{2n_1 - 1}) \ar[r] & 0, \\
0 \ar[r] & \K_2 \ar[r] & A_2 \ar[r]^-{\sigma_{\Rr^{n_2}}} & C(S^{2n_2 - 1}) \ar[r] & 0.
}
\end{align*}

Note that the index map in $K$-theory in each case is surjective, e.g. by Fedosov's index formula.
From the six-term exact sequence we see that the index map from $\Zz \to \Zz$ is in fact an isomorphism.

Next note that $A_1, A_2$ are separable and as extensions of nuclear $C^{\ast}$-algebras themselves nuclear.
Also as we have just seen the $K$-theory groups are torsion-free. 
Hence we can apply K\"unneth's theorem (\cite{wegge}, p. 171) as follows for $i = 0, \ 1$ and $j = 1, \ 2$:
\begin{align*}
& K_i(C(S^{2n_j - 1}) \otimes A_1) = K_i(C(S^{2n_j-1}) \otimes A_2) \cong \Zz, \\ 
& K_0(A_1 \otimes \K_2) = K_0(A_2 \otimes \K_1) \cong \Zz, \\
& K_1(A_1 \otimes \K_2) = K_1(A_2 \otimes \K_1) \cong 0.
\end{align*}

It follows with Prop. \ref{Lem:1}, \emph{ii)}:
\[
K_0(A^{-1,0}) = K_0(A^{0,-1}) \cong \Zz, \ K_1(A^{0,-1}) = K_1(A^{-1,0}) \cong 0.
\]

\emph{ii)} We calculate the $K$-theory of the pullback $\Sigma$ via Mayer-Vietoris in $K$-theory (cf. \cite{wegge}, 11.D).

This gives:
\[
\xymatrixrowsep{0.5in}
\xymatrixcolsep{0.5in}
\xymatrix{
K_0(\Sigma) \ar[r]^-{\pi_{1^{\ast}} \oplus \pi_{2^{\ast}}} & K_0(-, A_1) \oplus K_0(-, A_2) \ar[r]^-{\tsigma_{{\Rr^{n_1}}^{\ast}} - \tsigma_{{\Rr^{n_2}}^{\ast}}} & K_0(C(S^{2n_1-1} \times S^{2n_2 -1})) \ar[d]_{\epsilon} \\
K_1(C(S^{2n_1-1} \times S^{2n_2 -1})) \ar[u]^{\delta} & \ar[l]^-{\tsigma_{\Rr^{n_1}}^{\ast} - \tsigma_{\Rr^{n_2}}^{\ast}} K_1(-, A_1) \oplus K_1(-, A_2) & \ar[l]^-{\pi_1^{\ast} \oplus \pi_2^{\ast}} K_1(\Sigma) 
}
\]

Here $K_i(-, A_1) = K_i(C(S^{2n_2 - 1}, A_1)), \ K_i(-, A_2) = K_i(C(S^{2n_1 - 1}, A_2), i = 0,1$.

With the isomorphisms already established in \emph{i)} we just have to calculate the maps on the generators of $\Zz \oplus \Zz$ 
in each case.
So we denote by $[1]_0, \ [\tilde{1}]_0$ the generators of $K^0(S^{2n_1 - 1}), \ K^0(S^{2n_2 - 1})$ respectively.
As well as by $[u]_1, \ [\tilde{u}]_1$ the unitary generators of $K^1(S^{2n_1 - 1}), \ K^1(S^{2n_2 - 1})$ respectively.

Then $K_0(\tsigma_{\Rr^{n_1}}) - K_0(\tsigma_{\Rr^{n_2}}) = \tsigma_{{\Rr^{n_1}}^{\ast}} - \tsigma_{{\Rr^{n_2}}^{\ast}}$ is given by
\begin{align}
&\Zz \oplus \Zz \ni (k_0,l_0) \mapsto (k_0 - l_0, 0). \label{first}
\end{align}

This follows since from the short exact sequence $0 \to \K_1 \to A_1 \to C(S^{2n_1 - 1}) \to 0$ and similarly for $A_2$, we 
see that the generators of $K_0(A_1), \ K_0(A_2)$ are determined by $[1]_0, \ [\tilde{1}]_0$ respectively.

Secondly, the map $K_1(\tsigma_{\Rr^{n_1}}) - K_1(\tsigma_{\Rr^{n_2}}) = \tsigma_{\Rr^{n_1}}^{\ast} - \tsigma_{\Rr^{n_2}}^{\ast}$ is given by
\begin{align}
\Zz \oplus \Zz \ni (k_1, l_1) \mapsto (k_1, -l_1) \in \Zz \oplus \Zz. \label{second}
\end{align}

To prove this observe that 
\begin{align*}
K_1(C(S^{2n_1 - 1}) \otimes A_2) &\cong K_1(C(S^{2n_1 - 1})) \otimes K_0(A_2) \cong \Zz, \\
K_1(C(S^{2n_2 - 1}) \otimes A_1) &\cong K_1(C(S^{2n_2 - 1})) \otimes K_0(A_1) \cong \Zz.
\end{align*}

Then $\tsigma_{\Rr^{n_1}}^{\ast}$ maps $[u]_1 \otimes [\tilde{1}]_0$ to $([u]_1 \otimes [\tilde{1}]_0, 0)$ and $\tsigma_{\Rr^{n_2}}^{\ast}$ maps $[\tilde{u}]_1 \otimes [1]_0$ to $(0, [\tilde{u}]_1 \otimes [1]_0)$. 

The morphism \eqref{second} is an isomorphism and hence the preceding and following arrows are zero maps. 
With this we have:

\[
\xymatrixrowsep{0.4in}
\xymatrixcolsep{0.8in}
\xymatrix{
K_0(\Sigma) \ar[r]^{} & \Zz \oplus \Zz \ar[r]^{(k_0,l_0) \mapsto (k_0 - l_0, 0)} & \Zz \oplus \Zz \ar[d]_{\epsilon} \\
\Zz \oplus \Zz \ar[u]^{0} & \ar[l]^{(k_1, -l_1) \mapsfrom (k_1,l_1)} \Zz \oplus \Zz & \ar[l]^{0} K_1(\Sigma) 
}
\]

It follows
\begin{align*}
K_0(\Sigma) &\cong \ker\{\Zz \oplus \Zz \ni (k_0,l_0) \mapsto (k_0-l_0, 0) \in \Zz \oplus \Zz\} \cong \Delta \cong \Zz, \\
K_1(\Sigma) &\cong \coker\{\Zz \oplus \Zz \ni (k_0,l_0) \mapsto (k_0-l_0,0) \in \Zz \oplus \Zz\} \cong \Zz^2 / \Zz \cong \Zz,
\end{align*}

where $\Delta$ denotes the diagonal in $\Zz^2$.

Finally, we determine the $K$-theory of $\A \cong A_1 \otimes A_2$ by again applying K\"unneth's theorem.
This yields $K_0(\A) \cong \Zz, \ K_1(\A) \cong 0$. 
\end{proof}

% prove that there is index 1 element in Shubin calculus

\begin{Rem}
There is an alternative way to obtain the $K$-theory of the symbol algebra; relying on the existence of an element of index one 
in the bisingular algebra.

First, we recall the following well-known property: Let $A \subset \L(\H)$ be a $C^{\ast}$-algebra containing the compact operators $\K$ and containing an element $a \in A$
with Fredholm index one. 
Then it is a well-known fact that the $K$-theory of the quotient $A / \K$ is given by 
\[
K_0(A) \cong K_0(A / \K), \ K_1(A) \oplus \Zz \cong K_1(A / \K). 
\]

To see this apply the six-term exact sequence in $K$-theory to the short exact sequence 
\[
\xymatrix{
\K \ar[r]^{j} & A \ar[r]^-{q} & A / \K.
}
\]

Since $A$ contains the compacts and an element of index one there is a non-unitary isometry which acts with regard to some fixed
orthonormal basis.
Therefore every finite rank projection in the compacts is stably homotopic to $0$ and the map in $K$-theory induced by the inclusion $j$ is the zero map.
This proves the claim. 

Applying this result to our bisingular algebra $\A$ we first note that with $\A \cong A_1 \otimes A_2$ and by the K\"unneth theorem $K_0(\A) \cong K_0(A_1) \otimes K_0(A_2) \oplus K_1(A_1) \otimes K_1(A_2)
\cong \Zz$ as well as $K_1(A) \cong K_1(A_1) \otimes K_0(A_2) \oplus K_0(A_1) \otimes K_1(A_2) \cong 0$.  
The Shubin classes contain elements $P_1 \in A_1$ and $P_2 \in A_2$ each of Fredholm index one, see e.g. \cite{hoerm}, Thm. 19.3.1. 
The external product $P_1 \sharp P_2$ is contained in the bisingular calculus $\A$ and the multiplicativity of the Fredholm index yields $\ind(P_1 \sharp P_2) = \ind(P_1) \cdot \ind(P_2) = 1$.
Hence by the above result we obtain the $K$-theory $K_0(\A / \K) \cong \Zz, \ K_1(\A/ \K) \cong K_0(\K) \cong \Zz$. 
\label{Rem:Kthy}
\end{Rem}

% mention general result: existence operators of order 1 and contains compact operators => K-thy \Zz and \Zz \oplus K_1(\A)
% std argument, non-unitary isometries 

\section{Toeplitz operators}

% Formulate general index problem
% Prop.: Index theory same for several classes of Toeplitz operators
% Special case from K-theory section: Globally MV sequence
% local vs. global index
% deduce abstract nonsense index thm

The index problem for the completed bisingular operator algebra can be stated 
abstractly without recourse to any pseudodifferential calculus.

At first we will introduce the classes of \emph{admissible Toeplitz operators} and calculate their $K$-theory. 
It also holds that the index theory of tensor products of such operators is the same as the index theory for the class of 
global bisingular operators.

% introduce admissible toeplitz classes

% general scheme: \Tau_{\H} = P_{\H} \pi(A) P_{\H} + \K(\H) 

\textbf{\emph{General scheme:}} we are given a Hilbert space $\tilde{\H}$ and a closed subspace $\H \subset \tilde{\H}$. 
Futhermore, $A$ is a $C^{\ast}$-algebra and $\pi \colon A \to \L(\H)$ is a representation of $A$. 
We fix the orthogonal projection $P_{\H} \colon \tilde{\H} \to \H$ and define the space of general \emph{Toeplitz operators}
\begin{align}
\Tau_{\H} := P_{\H} \pi(A) P_{\H} + \K(\H). \label{Toeplitz}
\end{align}

In the following we consider three particular cases which we summarize as follows.
For further details the reader can consult the references \cite{folland} and \cite{shubin}. 

\begin{itemize}

\item \emph{Hardy:} $\tilde{\H}_S = L^2(S^{2n -1 }), \ \H_S = H^2(S^{2n-1}), \ A_S = C(\overline{B}_{2n})$. The orthogonal projection
is defined as
\[
(P_S f)(z) := \int_{S^{2n-1}} \frac{f(w)}{(1 - z \overline{w})^n} \,dS(w). 
\]

\item \emph{Bargmann:} $\tilde{\H}_B = C(B_{2n}), \ \H_B = B^2(B_{2n}), \ A_B = C(\overline{B}_{2n})$. The orthogonal projection is defined as
\[
(P_B f)(z) := \int_{B_{2n}} \frac{f(w)}{(1 - z \overline{w})^{n+1}}\,dV(w). 
\]

\item \emph{Fock:} $\tilde{\H}_{\F} = L^2(\Cc^n, e^{-|z|^2}), \ \H_{\F} = \F_n = B^2(\Cc^n, e^{-|z|^2}), \ A_{\F} = C(\overline{B}_{2n})$. 
The orthogonal projection is defined as
\[
(P_{\F} f)(z) := \int_{\Cc^n} f(w) e^{z \overline{w} - |w|^2}\,dV(w). 
\]

\end{itemize}

\begin{Def}
Let $I \in \{B, \F, S\}$, then 
\[
\Tau_I := P_{I} \pi(A_{I}) P_{I} + \K(\H_I)
\]

is called algebra of \emph{admissible Toeplitz operators}. 

\label{Def:admToeplitz}
\end{Def}

\begin{Thm}
Let $\Tau_1, \ \Tau_2$ be two algebras belonging to one class of admissible Toeplitz operators over $\Rr^{n_1}$ and $\Rr^{n_2}$ respectively.
Then the index theory of the $C^{\ast}$-tensor product $\Tau_1 \otimes \Tau_2$ is the same as the index theory of the
completed global bisingular operators.
\label{Thm:Toeplitz}
\end{Thm}

\begin{proof}
Let $I \in \{B, \F, S\}$ and denote by $\Tau_1, \ \Tau_1$ two algebras of admissible Toeplitz operators depending on $I$.
Fix two unitary operators $U_1 \colon \H_{1,I} \to L^2(\Rr^{n_1})$ and $U_2 \colon \H_{2,I} \to L^2(\Rr^{n_2})$
and denote by $\Ad(U_i) \colon L^2(\Rr^{n_i}) \to \L(\H_{i,I}), \ i = 1,2$ the induced isomorphism of algebras given by
\[
\Ad(U_i) \colon \L(L^2(\Rr^{n_i})) \to \L(\H_{i,I}), \ i = 1,2, \ T \mapsto U^{\ast} T U. 
\]

We denote by $\K_i(I) = \K(\H_{i,I})$ the compact operators on the Hilbert spaces $\H_{i,I}$ for $i = 1,2$.

For $i = 1,2$ we have the commuting diagram (see \cite{folland})
\[
\xymatrix{
\K(L^2(\Rr^{n_i})) \ar[d]_{\Ad(U_i)} \ar@{>->}[r] & \overline{G_{cl}^{0}(\Rr^{n_i})} \ar[d]_{\Ad(U_i)} \ar@{->>}[r] & C(S^{2n_i - 1}) \ar@{--}[d] \\
\K_i(I) \ar@{>->}[r] & \Tau_i \ar@{->>}[r] & C(S^{2n_i - 1})
}
\]

We obtain the commuting diagram
\[
\xymatrix{
\K_1(I) \otimes \K_2(I) \ar[d]_{\Ad(U_1) \otimes \Ad(U_2)} \ar@{>->}[r] & \Tau_1 \otimes \Tau_2 \ar[d]_{\Ad(U_1) \otimes \Ad(U_2)} \ar@{->>}[r] & \Sigma_{\Tau} \ar[d]_{\Phi} \\
\K_1 \otimes \K_2 \ar@{>->}[r] & \overline{G_{cl}^{0}(\Rr^{n_1})} \otimes \overline{G_{cl}^{0}(\Rr^{n_2})} \ar@{->>}[r] & \Sigma
}
\]

Let's describe the $\ast$-homomorphism $\Phi$. For this first write down the two pullbacks
\[
\xymatrix{
\Sigma_{\Tau} \ar[d]_{\pi_2} \ar[r]^-{\pi_1} & C(S^{2n_1 - 1}, \Tau_2) \ar[d]_{\id_1 \otimes \tilde{\sigma}_2^{I}} \\
C(S^{2n_2 - 1}, \Tau_1) \ar[r]^-{\id_2 \otimes \tilde{\sigma}_1^{I}} & C(S^{2n_1 - 1} \times S^{2n_2 - 1})
}
\]

and
\[
\xymatrix{
\Sigma \ar[d]_{\pi_2} \ar[r]^-{\pi_1} & C(S^{2n_1 - 1}, \overline{G_{cl}^{0}(\Rr^{n_2})}) \ar[d]_{\id_1 \otimes \sigma_2} \\
C(S^{2n_2 - 1}, \overline{G_{cl}^{0}(\Rr^{n_1})}) \ar[r]^{\id_2 \otimes \sigma_1} & C(S^{2n_1 - 1} \times S^{2n_2 - 1}). 
}
\]

Then we define the map $\Phi \colon \Sigma_{\Tau} \to \Sigma$ by 
\[
\Phi := (\id_1 \otimes \Ad(U_2)) \oplus (\Ad(U_1) \otimes \id_2). 
\]

We check that $\Phi$ is a well-defined $\ast$-homomorphism.
To this end consider the following diagram
\SelectTips{eu}{12}
\[ \xymatrix{ C(S^{2n_2 - 1}, \Tau_1) \ar[dd]_{\id_2 \otimes \Ad(U_1)} \ar[rr]^{\id_2 \otimes \tilde{\sigma}_1^{I}} && C(S^{2n_1 - 1} \times S^{2n_2 - 1}) \ar@{==}'[d][dd] \\
& \ar[ul]^{\tilde{\pi}_2} \Sigma_{\Tau} \ar@{-->}[dd] \ar[rr]^{\tilde{\pi}_1} && \ar[ul]^{\id_1 \otimes \tilde{\sigma}_2^I} C(S^{2n_1 -1}, \Tau_2) \ar[dd]_{\id_1 \otimes \Ad(U_2)} \\
C(S^{2n_2 -1}, A_1) \ar'[r][rr]^-{\id_2 \otimes \tilde{\sigma}_1} && C(S^{2n_1 -1} \times S^{2n_2-1}) \\
& \ar[ul]_{\pi_2} \Sigma \ar[rr]^{\pi_1} && \ar[ul]^{\id_1 \otimes \tilde{\sigma}_2} C(S^{2n_1 - 1}, A_2) }
\]

Let $(F, G) \in \Sigma_{\Tau}$, then by definition
\[
(\id_1 \otimes \tilde{\sigma}_2^{I})(F) = (\tilde{\sigma}_1^{I} \otimes \id_2)(G). 
\]

For the following we can wlog set $F = f \otimes a$ for $f \in C(S^{2n_1 - 1}), \ a \in \Tau_2$ and $G = b \otimes g, \ g \in C(S^{2n_2 -1}), \ b \in \Tau_1$. 

Then we check that
\begin{align*}
&(\id_1 \otimes \sigma_2)(\id_1 \otimes \Ad(U_2))(F) \\
&= (\id_1 \otimes \sigma_2)(f \otimes U_2^{\ast} a U_2) \\
&= f \otimes \sigma_2(U_2^{\ast} a U_2) \\
&= f \otimes \tilde{\sigma}_2^{I}(a) \\
&= \tilde{\sigma}_1^{I}(b) \otimes g \\
&= (\sigma_1 \otimes \id_2)(U_1^{\ast} b U_1 \otimes g) \\
&= (\sigma_1 \otimes \id_2)(\Ad(U_1) \otimes \id_2)(G). 
\end{align*}

Hence $\Phi(F, G) \in \Sigma$. 
The commutativity of the diagram follows from a similar calculation.

This suffices to prove the assertion concerning the index theory of these operators.
Denote by
\[
\partial_I \colon K_1(\Sigma_{\Tau}) \to \Zz = K_0(\K)
\]

the connecting map in $K$-theory. 
Then for $T \in \Tau_1(I) \otimes \Tau_2(I) \otimes M_N(\Cc)$ a given Fredholm operator we have that
\[
\ind(T) = \partial_I[(\sigma_1^I \oplus \sigma_2^I)(T)]_1
\]

and by the previous result we can calculate the index independent of the choice of $I$. 
% not quite: define \Phi 
% show: \Phi well-defined \ast-isomorphism
% commutativity of the diagram?
\end{proof}

The $K$-theory of tensor products of the Toeplitz operators can be calculated completely analogously by use of the
K\"unneth theorem and the Mayer-Vietoris sequence. 
\newpage

% corollary: table with K-theory
\begin{table}[ht]
\caption{$K$-theory}
\centering
\begin{tabular}{c c c}
\hline \hline
$C^{\ast}$-algebra & $K_0$ & $K_1$ \\ [0.5ex]
\hline
$\Tau_1 \otimes \Tau_2$ & $\Zz$ & $0$ \\
$\Tau_1 \otimes \Tau_2 / \K_1(I) \otimes \K_2(I)$ & $\Zz$ & $\Zz$ \\ 
$\overline{G_{cl}^{-1,-1}} = \K$ & $\Zz$ & $0$ \\
$\overline{G_{cl}^{-1,0}}$ & $\Zz$ & $0$ \\
$\overline{G_{cl}^{0,-1}}$ & $\Zz$ & $0$ \\
$\overline{G_{cl}^{0,0}} / \K$ & $\Zz$ & $\Zz$ \\ [1ex]
\hline
\end{tabular}
\label{table:Kthy}
\end{table}

% abstract nonsense approach using KK-theory

% A_1, A_2 .. \tilde{\sigma}_1, \tilde{\sigma}_2 

\section{Concluding remarks}

\emph{i)} We consider the general setup for index theory as in section \ref{Kthy}.

% A_1, A_2 .. \tilde{\sigma}_1, \tilde{\sigma}_2 

Denote by $\E_A, \E_B$ two separable $C^{\ast}$-algebras and by $A, B$ nuclear, separable $C^{\ast}$-algebras. 
Such that we have the two exact sequences (with completed $C^{\ast}$-tensor products)
\begin{align*}
\xymatrix{
\K_1 \ar@{>->}[r] & \E_A \ar@{->>}[r]^{\sigma_A} & A, \\
\K_2 \ar@{>->}[r] & \E_B \ar@{->>}[r]^{\sigma_B} & B 
}
\end{align*}

and we have the $C^{\ast}$-pullback 
\begin{align*}
\xymatrix{
\Sigma \ar[d]_-{\id_A \otimes \sigma_B} \ar[r]^-{\sigma_A \otimes \id} & A \otimes \E_B \ar[d]_{\id_A \otimes \sigma_B} \\
\E_A \otimes B \ar[r]^{\sigma_A \otimes \id_B} & A \otimes B
}
\end{align*}

Then from the two exact sequences it follows that $\E_A, \E_B$ are nuclear as extensions of nuclear algebras.
Additionally, $\Sigma$ fits into two exact sequences
\[
\xymatrix{
\K_1 \otimes B \ar@{>->}[r] & \Sigma \ar@{->>}[r] & A \otimes \E_B 
}
\]

and
\[
\xymatrix{
\K_2 \otimes A \ar@{>->}[r] & \Sigma \ar@{->>}[r] & B \otimes \E_A.
}
\]

So in particular $\Sigma$ is nuclear. 

On can then consider the generalized index problem: \emph{Express the generalized analytic index of the 
tensor product algebra $\E_A \otimes \E_B$ in topological terms}.

After application of the Mayer-Vietoris sequence in $K$-theory we consider the following diagram
\begin{align*}
\xymatrix{
& \ar[dl]_{\ind_a} K_1(\Sigma) \\
\Zz & \\
& \ar[ul]_{\ind_{loc}} K_0(A \otimes B) \ar[uu]_{\partial_{MV}} & 
}
\end{align*}

Here $\partial_{MV}$ denotes the connecting map in Mayer-Vietoris as defined in Theorem \ref{Thm:MV}, $\ind_a$ is given 
by the connecting map in $K$-theory and $\ind_{loc}$ is the \emph{local index} which we are going to define next.

We set
\begin{align}
\ind_{loc} := - \otimes ([\E_A] \otimes [\E_B]) \label{loc}
\end{align}

where $[\E_A] \in \Ext^{-1}(\E_A) \cong KK_1(A, \Cc) \cong K^1(A)$ and $[\E_B] \in \Ext^{-1}(\E_B) \cong KK_1(B, \Cc) \cong K^1(B)$ denote the corresponding extension classes.
Here $\otimes$ denotes in each case the corresponding and appropriate Kasparov-product. 

Note that by \ref{Thm:MV}, \emph{ii)} the Mayer-Vietoris connecting map is given by $\partial_{MV} = - \otimes \partial_{\Sigma}$ for a $\partial_{\Sigma} \in KK(A \otimes B, \Sigma)$.

\emph{ii)} Consider the index for the global bisingular calculus (and therefore by Theorem \ref{Thm:Toeplitz} for tensor products
of admissible Toeplitz operators).

With the $K$-theory calculation of \ref{Thm:Kthy} and abstract nonsense we can obtain an expression of the analytical index.
In a work in preparation we make this expression explicit in terms of a Fedosov type index formula. 

There is a splitting $\beta \colon K_1(\Sigma) \to K^0(S^{2n_1 - 1} \times S^{2n_2 - 1})$ such that $\epsilon \circ \beta = \id$.
Then the topological index defined by
\[
\ind_t(a) := \int_{S^{2n_1 - 1} \times S^{2n_2 - 1}} \ch([\beta(a)]_0)
\]

agrees with the analytical index. 

This follows from the Mayer Vietoris sequence in the proof of \ref{Thm:Kthy}. 
Set for example $\beta(l) = (l \cdot m, l)$ for some fixed $m \in \Zz$. Since $\epsilon$ is surjective with constant kernel
we have $\epsilon(k,l) = l$ for $k,l \in \Zz$.
Then fix a positive orientation on $S^{2n_1 - 1} \times S^{2n_2 - 1}$ to define the topological index.
The equality follows then by exhibiting an operator in the bisingular calculus of index $1$.
For the latter use a standard construction or simply take the external product $P := P_1 \sharp P_2$ of two
index $1$ Shubin operators $P_1, P_2$. Then by multiplicativity $\ind(P_1 \sharp P_2) = \ind(P_1) \cdot \ind(P_2) = 1$.

\begin{appendix}

\section{The Mayer-Vietoris Theorem}

% appendix: Mayer Vietoris sequence
% \cite{black}, p. 219, section 21.2, Thm. 21.2.2

In this section we will state and prove the Mayer-Vietoris Theorem in $C^{\ast}$-algebra $K$-theory. 
The following proof can be found e.g. in \cite{black} in section 21.2. 
We nevertheless give some details for the benefit of the non-specialist and in order to make the paper more self-contained.

First we recall some standard constructions.

\begin{Def}
\emph{i)} Let $A$ be a $C^{\ast}$-algebra. Then we define
\begin{align*}
SA &:= \{f \in C([0,1], A) |f(0) = f(1) = 0\}, \\
\C A &:= \{f \in C([0,1], A) | f(0) = 0\}
\end{align*}

which is called the \emph{suspension} and \emph{cone} of $A$ respectively. 

\emph{ii)} Let $A, B$ be two $C^{\ast}$-algebras and $\varphi \colon A \to B$ a $\ast$-homomorphism. 
Then
\begin{align*}
\C_{\varphi} &:= \{(a, f) \in A \oplus \C B : \varphi(a) = f(1)\} 
\end{align*}

is called the \emph{mapping cone} of $\varphi$. 
\label{Def:SA}
\end{Def}

\begin{Rem}
We can immediately write down some standard exact sequences:
\begin{align}
\xymatrix{
0 \ar[r] & S B \ar[r]^-{i} & \C_{\varphi} \ar[r]^-{j} & A \ar[r] & 0 
} \label{SA1} 
\end{align}

with $i(f) := (0, f), \ j(a,f) = a$. 

Also
\begin{align}
\xymatrix{
0 \ar[r] & SA \ar[r] & \C A \ar[r]^{\pi} & A \ar[r] & 0
} \label{SA2}
\end{align}

with $\pi(f) := f(1)$ the evaluation. 

\label{Rem:SA}
\end{Rem}

\begin{Def}
Two $\ast$-homomorphisms $\varphi, \ \psi \colon A \to B$ are said to be homotopy-equivalent if there is a family
of $\ast$-homomorphisms $\omega_t \colon A \to B, \ t \in [0,1]$ with $t \mapsto \omega_t$ pointwise norm-continuous and
$\omega_0 = \varphi, \ \omega_1 = \psi$. We write $\varphi \sim_h \psi$. 
\label{Def:homotopy}
\end{Def}

\begin{Rem}
We have that $\varphi \sim_h \psi$ is equivalent to $\exists \ \omega \colon A \to C([0,1], B)$ with
\[
\delta_0 \circ \omega = \varphi, \ \delta_1 \circ \omega = \psi 
\]

where $\delta_t \colon C([0,1], B) \to B$ is evaluation at $t$. 
\label{Rem:homotopy}
\end{Rem}

% A \sim_h B :<=> ...

\begin{Prop}
The cone $\C A$ is contractible for each $C^{\ast}$-algebra $A$.
\label{Prop:CA}
\end{Prop}

\begin{proof}
Define $f_t \colon \C A \to \C A$ by $f_t(f)(s) := f(ts)$. Then $f_0 = 0, \ f_1 = \id_{CA}$, hence $\C A$ is contractible. 
\end{proof}

\begin{Thm}[Mayer-Vietoris]
Let $\E_1, \ \E_2$ and $\F$ be $C^{\ast}$-algebras and let $\varphi_1 \colon \E_1 \to \F, \ \varphi_2 \colon \E_2 \to \F$ be surjective $\ast$-homomorphisms. 
Denote by $\Sigma$ the $C^{\ast}$-pullback 
\[
\xymatrix{
\Sigma \ar[d]_-{\pi_2} \ar[r]^-{\pi_1} & \E_1 \ar[d]_-{\varphi_1} \\
\E_2 \ar[r]^-{\varphi_2} & \F
}
\]

which is the restricted direct sum $\Sigma = \E_1 \oplus_{\F} \E_2 := \{(a, b) \in \E_1 \oplus \E_2 : \varphi_1(a) = \varphi_2(b)\}$. 

\emph{i)} There is a six-term exact sequence in $K$-theory
\[
\xymatrix{
K_0(\Sigma) \ar[r]^-{(\pi_1 \oplus \pi_2)_{\ast}} & K_0(\E_1) \oplus K_0(\E_2) \ar[r]^-{\varphi_{1\ast} - \varphi_{2\ast}} & K_0(\F) \ar[d]_{\epsilon} \\ 
K_1(\F) \ar[u]^{\delta} & \ar[l]^-{\varphi_{1\ast} - \varphi_{2\ast}} K_1(\E_1) \oplus K_1(\E_2) & \ar[l]^-{(\pi_1 \oplus \pi_2)_{\ast}} K_1(\Sigma).
}
\]

\emph{ii)} The connecting map $\epsilon$ takes the form of the Kasparov product $- \otimes \partial_{\Sigma}$ for some
$\partial_{\Sigma} \in KK_1(\F, \Sigma)$. 
\label{Thm:MV}
\end{Thm}

\begin{proof}
\emph{i)} Let $\pi := \pi_1 \oplus \pi_2 \colon \Sigma \to \E_1 \oplus \E_2$, then the mapping cone
\[
\C_{\pi} = \{(a, f) \in \Sigma \oplus \C(\E_1 \oplus \E_2) : \pi(a) = f(0)\}
\]

is $\ast$-isomorphic to the $C^{\ast}$-algebra
\[
\{(h_1, h_2) | \varphi_1(h_1(0)) = \varphi_2(h_2(0))\} \subset \C \E_1 \oplus \C \E_2.
\] 

Define $\psi \colon \C_{\pi} \to S \F$ via 
\[
\psi(h_1, h_2)(t) := \begin{cases} \varphi_1(h_1(1 - 2t)), &\text{for} \ 0 \leq t \leq \frac{1}{2} \\
\varphi_2(h_2(2t - 1)), &\text{for} \ \frac{1}{2} \leq t \leq 1.
\end{cases}
\]

Then $\psi$ is a well-defined $\ast$-homomorphism. We can show that $\psi$ is also surjective as follows.
Let $f \in S \F$, then by surjectivity of $\varphi_1, \ \varphi_2$ choose $a \in \E_1, \ b \in \E_2$ such that 
$f(\frac{1}{2}) = \varphi_1(a) = \varphi_2(b)$. 
Therefore $(a, b) \in \Sigma$ and setting $h_1(t) = (1 - t) a, \ h_2(t) = (1 -t) b, \ t \in [0,1]$. 
Then $(h_1, h_2) \in \C_{\pi}$ and $\psi(h_1, h_2) = f$, therefore $\psi$ is surjective. 

Setting $J_i := \ker \varphi_i, \ i = 1,2$ it follows that
\[
\ker \psi = \C J_1 \oplus \C J_2. 
\]

We obtain the short exact sequence
\[
\xymatrix{
\C J_1 \oplus \C J_2 \ar@{>->}[r] & \C_{\pi} \ar[r]^-{\psi} & S \F. 
}
\]

Applying the six-term exact sequence in $K$-theory to this yields $K_j(\C_{\pi}) = K_j(S\F)$ for $j = 1,2$ by
the contractibility of $\ker \psi$, Proposition \ref{Prop:CA} and homotopy invariance of the $K$-functor. 

Consider the standard exact sequence (see \ref{Rem:SA}, \eqref{SA1})
\[
\xymatrix{
S \E_1 \oplus S \E_2 \ar@{>->}[r]^-{j} & \C_{\pi} \ar@{->>}[r]^-{\tilde{\psi}} & \Sigma.
}
\]

Applying the six-term exact sequence we obtain
\[
\xymatrix{
K_1(\E_1) \oplus K_1(\E_2) \ar[r]^-{r_{\ast}} & K_0(\C_{\pi}) \ar[r]^-{\tilde{\psi}_{\ast}} & K_0(\Sigma) \ar[d]_-{(\pi_1 \oplus \pi_2)_{\ast}} \\
K_1(\Sigma) \ar[u]^-{(\pi_1 \oplus \pi_2)_{\ast}} & \ar[l]^-{\tilde{\psi}_{\ast}} K_1(\C_{\pi}) & \ar[l]^-{r_{\ast}} K_0(\E_1) \oplus K_0(\E_2)
}
\]

where we note that we have the natural isomorphisms $K_j(S \E_i) \cong K_{j+1}(\E_i)$ for $i = 1,2$ and $j = 1,2$ by Bott periodicity. 
Noting also that $K_j(\C_{\pi}) \cong K_{j}(S\F) \cong K_{j+1}(\F)$ we obtain the desired form of the Mayer-Vietoris sequence.

Finally, the induced map in $K$-theory $r_{\ast}$ is given as follows. Setting $r := \psi_{|S \E_1 \oplus S \E_2}$ we find
that the function $\omega \colon S \E_1 \oplus S \E_2 \to S \F$ given by
\[
\omega(h_1, h_2)(t) := \varphi_1(h_2(t)) - \varphi_2(h_1(t))
\]

is homotopic to $r$. In particular $r_{\ast}(x, y) = \varphi_{2 \ast}(x) - \varphi_{1 \ast}(y)$. 

\emph{ii)} Using the notation from part \emph{i)} we consider the diagram
\[
\xymatrix{
S \E_1 \oplus S \E_2 \ar@{>->}[r] & \C_{\pi} \ar@{==}[d] \ar@{->>}[r]^{\tilde{\psi}} & \Sigma \\
\C J_1 \oplus \C J_2 \ar@{>->}[r] & \C_{\pi} \ar@{->>}[r]^{\psi} & S \F \ar[u]^{\partial_{\Sigma}}
}
\]

where the map $\partial_{\Sigma}$ is the sought for element in $KK$-theory.

By contractibility of $\C J_1 \oplus \C J_2$ we obtain for each $C^{\ast}$-algebra $A$ an isomorphism in $KK$-theory
\[
\psi_{\ast} = - \otimes [\psi] \colon KK(A, \C_{\pi}) \iso KK(A, S \F). 
\]

Setting $A = \C_{\pi}$ yields an invertible element in $KK$-theory $[\psi] \in KK(\C_{\pi}, S \F)$. 
Hence there is a $[\psi]^{-1} \in KK(S \F, \C_{\pi})$ such that 
\[
[\psi]^{-1} \circ [\psi] = 1_{S \F}, \ [\psi] \circ [\psi]^{-1} = 1_{\C_{\pi}}. 
\]

Setting $\partial_{\Sigma} = [\psi]^{-1} \circ [\tilde{\psi}]$ we recover the connecting map in $K$-theory $\epsilon = - \otimes \partial_{\Sigma}$. 
\end{proof}

\end{appendix}

\end{document}